% ...
\documentclass{svmult}   
\usepackage{graphicx}    

\usepackage{amsfonts}
\usepackage{amsmath}

\begin{document}


\title*{Convergence of products of stochastic matrices with positive diagonals and the opinion dynamics background} 

\titlerunning{Products of stochastic matrices with positive diagonals}


\author{
Jan Lorenz
} 





\index{Lorenz J.}



\institute{Universit{\"a}t Bremen, Fachbereich Mathematik und
Informatik, Bibliothekstra{\ss}e, 28359 Bremen, Germany
\texttt{math@janlo.de}}



\newcommand{\n}{\underline n}
\newcommand{\N}{\mathbb N}
\newcommand{\R}{\mathbb R}
\newcommand{\Ci}{\mathcal{I}}
\newcommand{\Cj}{\mathcal{J}}
\smartqed

\maketitle
\renewcommand{\abstractname}{Abstract.}
\begin{abstract}
We present a convergence result for infinite products of
stochastic matrices with positive diagonals. We regard infinity of
the product to the left. Such a product converges partly to a
fixed matrix if the minimal positive entry of each matrix does not
converge too fast to zero and if either zero-entries are symmetric
in each matrix or the length of subproducts which reach the
maximal achievable connectivity is bounded.

Variations of this result have been achieved independently in
\cite{Lorenz2005}, \cite{Moreau2005} and \cite{Hendrickx2005}. We
present briefly the opinion dynamics context, discuss the
relations to infinite products where infinity is to the right
(inhomogeneous Markov processes) and present a small improvement
and sketch another.
\end{abstract}


\section{Introduction}

Consider $n$ persons that discuss an issue which can be
represented as a real number. Assume further that the persons
revise their opinions if they hear the opinions of others. Each
person finds his new opinion as a weighted arithmetic mean of the
opinions of others. This model of opinion dynamics has been
analyzed for the possibilities of consensus by DeGroot
\cite{DeGroot1974}. If these weights change over time we have an
inhomogeneous consensus process.

While the homogeneous process has strong similarities with a
homogeneous Markov chain, things get different when inhomogeneity
comes in. While a consensus process relies on row-stochastic
matrices multiplied from the left, a Markov process relies on
row-stochastic matrices multiplied from the right. And infinity to
the right is unfortunately not the same as infinity to the left if
we consider row-stochastic matrices. But nevertheless, both
processes fit in the common framework of infinite products of
row-stochastic matrices.

Consensus processes are only briefly touched in the context of
Markov chains \cite{Hartfiel1998}. Besides the early approaches of
opinion dynamics \cite{DeGroot1974, Lehrer1981} some results have
been made in the context of decentralized computation
\cite{Tsitsiklis1984}. Recently there have been some independent
works that study consensus processes and the underlying
matrix-products in the context of opinion dynamics
\cite{Krause2000, Hegselmann2002, Lorenz2003b}, multiagent systems
where agents try to coordinate \cite{Moreau2005} and flocking
where birds try to find agreement about their headings
\cite{Hendrickx2005}.

In this paper we want to analyze the structure that positive
diagonals deliver in inhomogeneous consensus and Markov processes.

\section{Consensus and Markov processes}

For $n\in\N$ we define $\n := \{1,\dots,n\}$.

Let $A(0), A(1), \dots$ be a sequence of row-stochastic matrices.
Unless otherwise stated we regard every matrix to be square and to
have dimension $n$.

For natural numbers $s<t$ we define a \emph{forward accumulation}
$A(s,t) = A(s)\dots A(t-1)$ and a \emph{backward accumulation}
$A(t,s) = A(t-1)\dots A(s)$. Thus $A(s,s+1)=A(s+1,s)=A(s)$ and
$A(s,s)$ is the identity.

We briefly explain the two paradigmatic example processes where
forward and backward accumulations play a role.

Let $p(0)$ be a stochastic row vector and $p_i(0)$ is the $i$th
component of p(0), which is equal to the proportion of probability
mass or population which is in state $i$ at the beginning. The
sequence of row vectors $p(t) := p(0)A(0,t)$ is thus an
\emph{inhomogeneous Markov process} and $A(t)_{[i,j]}$ determines
the \emph{transition} from state $i$ to $j$ at time step $t$. In
that context $A(t)$ is called \emph{transition matrix}.

Let $x(0)$ be a real column vector of opinions and $x_i(0)$ stands
for the initial opinion of person $i$. The sequence of vectors
$x(t) = A(t,0)x(0)$ is an \emph{inhomogeneous consensus process}
and $A(t)_{[i,j]}$ stands for the weight person $i$ gives to the
opinion of agent $j$ at time step $t$. In that context $A(t)$ is
called \emph{confidence matrix}.

To understand the convergence behavior of inhomogeneous Markov and
consensus processes the infinite products $A(0,\infty)$ and
$A(\infty,0)$ are of interest.

In this paper we focus on transitions and confidence matrices with
positive diagonals. Thus, we regard Markov processes where we have
always a positive probability to stay in one state and consensus
processes where persons have at least a little bit of
self-confidence.

In the next section we will see that the positive diagonal
together with the Gantmacher form of nonnegative matrices will
give us a good overview on the zero and positivity structure of
the processes.

In section \ref{sec:conv} we go on with a convergence theorem that
is built on this structure and conclude in section \ref{sec:disc}
with a small improvement and discussion on how to fulfill the
conditions of the theorem.

\section{The positive diagonal} \label{sec:posdiag}

We regard two nonnegative matrices $A,B$ to be of the same
\emph{type} $A \sim B$ if $a_{ij} > 0 \Leftrightarrow b_{ij}>0$.
Thus, if their zero-patterns are equal.

Let $A$ be a nonnegative matrix with a positive diagonal. For
indices $i,j\in\n$ we say that there is a \emph{path} $i
\rightarrow j$ if there is a sequence of indices
$i=i_1,\dots,i_k=j$ such that for all $l\in\underline{k-1}$ it
holds $a_{i_l,i_{l+1}}>0$. We say $i,j \in\n$ \emph{communicate}
if $i \rightarrow j$ and $j \rightarrow i$, thus $i
\leftrightarrow j$. In our case with positive diagonals there is
always a path from an index to itself, which we call
\emph{self-communicating} and thus $'\leftrightarrow'$ is an
equivalence relation. An index $i\in\n$ is called \emph{essential}
if for every $j\in\n$ with $i \rightarrow j$ it holds $j
\rightarrow i$. An index is called inessential if it is not
essential.

Obviously, $\n$ divides into disjoint self-communicating
equivalence classes of indices $\Ci_1,\dots,\Ci_p$. Thus, in one
class all indices communicate and do not communicate with other
indices. The terms essential and inessential thus extend naturally
to classes.

If we renumber indices with first counting the essential classes
and second the inessential classes with a class $\Ci$ before a
class $\Cj$ if $\Cj\rightarrow\Ci$ then we can bring every
row-stochastic matrix $A$ to the \emph{Gantmacher form}
\cite{Gantmacher1959}
\begin{equation} \label{gf}
\begin{bmatrix}
  A_1 &  &  & &  & 0 \\
   & \ddots &  &  &  &  \\
  0 &  & A_g &  &  &  \\
  A_{g+1,1} & \dots & A_{g+1,g} & A_{g+1} &  &  \\
  \vdots &  & \vdots & \vdots & \ddots &  \\
  A_{p,1} & \dots & A_{p,g} & A_{p,g+1} & \dots & A_p \\
\end{bmatrix}
\end{equation}
by simultaneous row and column permutations. The \emph{diagonal
Gantmacher blocks} $A_1,\dots,A_p$ in (\ref{gf}) are square and
irreducible. Irreducibility induces primitivity in the positive
diagonal case. For the \emph{nondiagonal Gantmacher blocks}
$A_{k,l}$ with $k=g+1,\dots,p$ and $l=1,\dots,k-1$ it holds that
for every $k\in\{g+1,\dots,p\}$ at least one block of
$A_{k,1},\dots,A_{k,k-1}$ contains at least one positive entry.

The following proposition shows that an infinite backward or
forward accumulation of nonnegative matrices can be divided after
a certain time step into subaccumulations with a common Gantmacher
form.

\begin{proposition} \label{prop1}
Let $(A(t))_{t\in\N_0}$ be a sequence of nonnegative matrices with
positive diagonals. Then for the backward accumulation there
exists a sequence of natural numbers $0 < t_0 < t_1 < \dots$ such
that for all $i\in\N_0$ it holds
\begin{equation}
A(t_{i+1},t_{i}) \sim A(t_1,t_0).
\end{equation}
Thus, $A(t_{i+1},t_{i})$ can be brought to the same Gantmacher
form for all $i\in\n$. Further on, all Gantmacher diagonal blocks
are positive and all nondiagonal Gantmacher-Blocks are either
positive or zero.

The same holds for another sequences $0 < s_0 < s_1 < \dots$ for
the forward accumulation and $A(s_i,s_{i+1})$.
\end{proposition}

\begin{proof}
(In sketch, for more details see \cite{Lorenz2005}.)

The proof works with a double monotonic argument on the positivity
of entries: While more and more (or exactly the same) positive
entries appear in $A(t,0)$ monotonously increasing with rising
$t$, we reach a maximum at $t^\ast_0$. We cut $A(t^\ast_0,0)$ of
and find $t^\ast_1$ when $A(t,t^\ast_0)$ reaches maximal
positivity again with rising $t$. We go on like this and get the
sequence $(A(t^\ast_{i+1},t^\ast_i))_{i\in\N_0}$. Obviously, less
and less (or exactly the same) positive entries appear
monotonously decreasing with rising $i$ and we reach a minimum at
$k$. We relabel $t_j := t^\ast_{k+j}$ and thus have the desired
sequence $(t_i)_{i\in\N_0}$ with $A(t_{i+1},t_i)$ having the same
zero-pattern.

Positivity of Gantmacher blocks follows for all blocks
$A(t_{i+1},t_i)_{[\Cj,\Ci]}$ where we have a path $\Cj \rightarrow
\Ci$. If we have such a path, then there is a path from each index
in $\Cj$ to each index in $\Ci$ and thus every entry must be
positive in a long enough accumulation. Thus, the block has to be
positive already, otherwise $(t_i)_{i\in\N_0}$ is chosen wrong.

To prove the result for forward accumulations, we can use the same
arguments. \qed
\end{proof}

The next section regards the convergence behavior of the
Gantmacher diagonal blocks.

\section{Convergence} \label{sec:conv}

A row-stochastic matrix $K$ which has rank 1 and thus equal rows
is called a \emph{consensus matrix} because for a real vector $x$
it holds that $Kx$ is a vector with equal entries and thus
represents consensus among persons in a consensus process. Suppose
that $A(t):=K$ is a consensus matrix. It is easy to see that for
all $u\geq t$ it holds for the backward accumulation that $A(u,0)
= K$ while for the forward accumulation it only holds that
$A(0,u)$ is a consensus matrix but may change with $u$.

We define the \emph{coefficient of ergodicity} of a row-stochastic
matrix $A$ according to Hartfiel \cite{Hartfiel1998} as
\[\tau(A) := 1 - \min_{i,j \in \n} \sum_{k=1}^n \min\{a_{ik}, a_{jk}\}.\]

The coefficient of ergodicity of a row-stochastic matrix can only
be zero, if all rows are equal, thus if it is a consensus matrix.

The coefficient of ergodicity is submultiplicative (see
\cite{Hartfiel1998}) for row-stochastic matrices $A_0,\dots,A_i$
\begin{equation}\label{eq:submulterg}
\tau(A_i \cdots A_1 A_0) \leq \tau(A_i) \cdots \tau(A_1)
\tau(A_0).
\end{equation}

If $\lim_{t\to\infty} \tau(A(0,t)) = 0$ we say that $A(0,t)$ is
\emph{weakly ergodic}. Weakly ergodic means that the $A(0,t)$ gets
closer and closer to the set of consensus matrices and thus the
Markov process gets totally independent of the initial
distribution $p(0)$.

For $M \subset \R_{\geq 0}$ we define $\min^{+}M$ as the smallest
positive element of $M$. For a stochastic matrix $A$  we define
$\min^{+} A := \min^{+}_{i,j \in\n} a_{ij}$. We call $\min^{+}$
the \emph{positive minimum}.

For the positive minimum of a set of row-stochastic matrices
$A_0,\dots,A_i$ it holds
\begin{equation}\label{eq:posmin}
\min\mbox{}^{+}(A_i\cdots A_0) \geq \min\mbox{}^{+}A_i\cdots
\min\mbox{}^{+}A_0.
\end{equation}

\begin{theorem} \label{theorem}
Let $(A(t))_{t\in\N_0}$ be a sequence of row-stochastic matrices
with positive diagonals, $0<t_0<t_1<\dots$ be the sequence of time
steps defined by proposition \ref{prop1}, $\Ci_1,\dots,\Ci_g$ be
the essential and $\Cj$ be the union of all inessential classes of
$A(t_1,t_0)$.

If for all $i\in \N_0$ it holds $\min^{+} A(t_{i+1},t_i) \geq
\delta_i$ and $\sum_{i=1}^\infty \delta_i = \infty$, then
\[
\lim_{t\to\infty} A(t,0) = \left[\begin{array}{ccc|c}
  K_1 &  & 0 & 0 \\
   & \ddots &  & \vdots \\
  0 &  & K_g & 0 \\ \hline
   & \mathrm{not\ converging} &  & 0 \\
\end{array}\right] A(t_0,0)
\]
where $K_1,\dots,K_g$ are consensus matrices. (The matrices have
to be sorted by simultaneous row and column permutations according
to $\Ci_1,\dots,\Ci_g,\Cj$.)
\end{theorem}

\begin{proof}
The interesting blocks are the diagonal blocks. It is easy to see
due to the lower block triangular Gantmacher form of
$A(t_{i+1},t_i)$ for all $i\in\N_0$, that all diagonal blocks only
interfere with themselves when matrices are multiplied.

Let us regard the essential class $\Ci_k$ and abbreviate $A_i :=
A(t_{i+1},t_i)_{[\Ci_k,\Ci_k]}$.

We show that the minimal entry in a column $j$ of a row-stochastic
matrix $B$ cannot sink when multiplied from the right with another
row-stochastic matrix $A$,
\[\min_{i\in\n} (AB)_{ij} = \min_{i\in\n} \sum_{k=1}^n a_{ik}b_{kj} \geq \min_{i\in\n} b_{ij}.\]
Thus, the minimum of entries in column $j$ of the product
$A_i\cdots A_0$ is monotonously increasing with rising $i\in\N_0$.
With similar arguments it follows that the maximum of entries in
column $j$ of the product $A_i\cdots A_0$ is monotonously
decreasing with rising $i\in\N_0$.

Further on, it holds due to (\ref{eq:submulterg}) and the
definition of the coefficient of ergodicity that
\[\lim_{i\to\infty}\tau(A_i\dots A_1 A_0) \leq
\prod_{i=1}^{\infty}\tau(A_i) = \prod_{i=1}^{\infty}(1-\delta_i)
\leq \prod_{i=1}^{\infty}e^{-\delta_i} =
e^{-\sum_{i=1}^{\infty}\delta_i} = 0.
\]
The maximal distance of rows shrinks to zero. Both arguments
together imply that $\lim_{i\to\infty}(A_i\dots A_1 A_0)$ is a
consensus matrix which we call $K_k$.

Now it remains to show that the $[\Cj,\Cj]$-diagonal block of the
inessential classes converges to zero.

Let us define $||\cdot||$ as the row-sum-norm for matrices. It
holds $||A_{[\Cj,\Cj]}(t_{i+1},t_{i})|| \leq (1-\delta_i)$ and
thus like above it holds
\[||A_{[\Cj,\Cj]}(\infty,t_0)|| \leq \prod_{i=1}^\infty ||A_{[\Cj,\Cj]}(t_{i+1},t_i)|| \leq \prod_{i=1}^\infty (1-\delta_i) \leq = 0.\]
This proves that $\lim_{t\to\infty}A_{[\Cj,\Cj]}(t,0) = 0$. \qed
\end{proof}

An inhomogeneous consensus process $A(t,0)x(0)$ with persons who
have some self-confidence stabilizes (under weak conditions) such
that we have $g$ consensual subgroups (the essential classes)
which have internal consensus, while all other persons (the
inessential indices) may hop still around building opinions as
convex combinations of the values reached in the consensual
groups.

We will not treat the Markov case in detail. But a similar result
can be made, but not with fixed matrices $K_1,\dots,K_g$ but with
weak ergodicity after the time step $s_0$ within the independent
subgroups $\Ci_1,\dots,\Ci_g$.

In an inhomogeneous Markov process $p(0)A(0,t)$ a certain number
of independent absorbing classes evolve which get independent of
their initial conditions after a certain time step.

\section{Discussion on conditions for $\min^{+}A(t_{i+1},t_i) \geq
\delta_i$}\label{sec:disc}

One thing where theorem \ref{theorem} stays unspecific is that it
demands lower bounds for the positive minimum of the accumulations
$A(t_{i+1},t_i)$. But, what properties of the single matrices may
ensure the assumption $\min^{+}A(t_{i+1},t_i) \geq \delta_i$ with
$\sum \delta_i = \infty$?

The first idea would be to assume a \emph{uniform lower bound for
the positive minimum} $\delta < \min^{+}A(t)$ for all $t$. But
this is not enough.

Recent independent research
\cite{Moreau2005,Hendrickx2005,Lorenz2003b} has shown that either
 \emph{bounded intercommunication intervals} ($t_{i+1}-t_i < N$
for all $i\in\N_0$) or \emph{type-symmetry} ($A \sim A^T$) of all
matrices $A(t)$ can be assumed additional to the uniform lower
bound for the positive minimum to ensure the assumptions of
theorem \ref{theorem}. But improvements are possible.

\subparagraph{Bounded intercommunication intervals} Let us regard
$\delta < \min^{+}A(t)$ for all $t\in\N_0$. If $t_{i+1}-t_i \leq
N$ it holds by (\ref{eq:posmin}) that $\min^{+}A(t_{i+1},t_i) \geq
\delta^{N}$ and thus $\sum_{i=0}^\infty \delta^N = \infty$ and
thus theorem \ref{theorem} holds. But $t_{i+1}-t_i$ may slightly
rise as the next two propositions show.

\begin{proposition}
Let $0<\delta<1$ and $a \in \R_{>0}$ then
\begin{equation}
\sum_{n=1}^\infty \delta^{a\log(n)} < \infty \Longleftrightarrow
\delta < e^{-1}.
\end{equation}
\end{proposition}

\begin{proof}
We can use the integral test for the series $\sum_{n=1}^\infty
\delta^{a\log(n)}$ because $f(x) := \delta^{a\log(x)}$ is positive
and monotonously decreasing on $[1,\infty[$.

With substitution $y = \log(x)$ (thus $dx=e^{y}dy$) it holds
\begin{eqnarray*}
\int_1^\infty \delta^{a\log(x)}dx & = & \int_1^\infty
e^{a\log(\delta)\log(x)}dx = \int_1^\infty
e^{a\log(\delta)y}e^{y}dy \\
& = & \int_1^\infty e^{ay(\log(\delta)+1)}dy
\end{eqnarray*}
The integral is finite if and only if $\log(\delta)+1 < 0$ and
thus if $\delta < e^{-1}$. \qed
\end{proof}

\begin{proposition}
Let $0<\delta<1$ and $a \in \R_{>0}$ then
\begin{equation}
\sum_{n=3}^\infty \delta^{a\log(\log(n))} = \infty.
\end{equation}
\end{proposition}

\begin{proof}
We can use the integral test for the series $\sum_{n=1}^\infty
\delta^{a\log(\log(n))}$ because $f(x) := \delta^{a\log(\log(x))}$
is positive and monotonously decreasing on $[3,\infty[$.

With substitution $y = \log(\log(x))$ (thus $dx=e^{(y+e^y)}dy$) it
holds
\begin{eqnarray*}
\int_3^\infty \delta^{a\log(\log(x))}dx & = & \int_1^\infty
e^{a\log(\delta)\log(\log(x))}dx = \int_1^\infty
e^{a\log(\delta)y}e^{y+e^y}dy \\
& = & \int_1^\infty e^{ay(\log(\delta)+1) + e^y}dy
\end{eqnarray*}
The integral diverges because  $ay(\log(\delta)+1) + e^y
\longrightarrow \infty$ as $y\to\infty$. \qed
\end{proof}

Thus, assuming $\min^{+}A(t) > \delta > 0$ for all $t\in\N_0$ we
can allow a slow growing of $t_{i+1} - t_i$ to fulfill the
assumptions of theorem \ref{theorem}. Acceptable is a growing as
quick as $\log(\log(i))$. If $t_{i+1} - t_i$ grows as $\log(i)$
then it must hold $\delta > e^{-1} > \frac{1}{3}$. This can only
hold if each row of $A(t)$ contains only two positive entries (due
to row-stochasticity).

\subparagraph{Type-symmetry}

Another way to ensure the assumptions of theorem \ref{theorem} is
to demand all matrices to by type-symmetric and have a positive
minimum uniformly bounded form below by $\delta$. Perhaps a small
improvement can be made with a very slowly sinking positive
minimum approaching zero. But giving a precise form is a task for
future work.

\subsection*{Acknowledgement}
I would like to thank Dirk Lorenz for calculus hints.

\bibliographystyle{unsrt}

%

\end{document}